\documentclass[12pt,leqno]{article}
\usepackage{amsmath,amssymb,amsthm,amscd,dsfont}
\numberwithin{equation}{section} \allowdisplaybreaks
\begin{document}
\newtheorem{theorem}{Theorem}[section]
\newtheorem{defin}{Definition}[section]
\newtheorem{prop}{Proposition}[section]
\newtheorem{corol}{Corollary}[section]
\newtheorem{lemma}{Lemma}[section]
\newtheorem{rem}{Remark}[section]
\newtheorem{example}{Example}[section]
\title{Weak-Hamiltonian dynamical systems}
\author{{\small by}\vspace{2mm}\\Izu Vaisman}
\date{}
\maketitle
{\def\thefootnote{*}\footnotetext[1]%
{{\it 2000 Mathematics Subject Classification: 53D20, 70H45, 93A30}.
}}
\begin{center} \begin{minipage}{12cm}
A{\footnotesize BSTRACT. A big-isotropic structure $E$ is an
isotropic subbundle of $TM\oplus T^*M$, endowed with the metric
defined by pairing. The structure $E$ is said to be integrable if
the Courant bracket $[\mathcal{X},\mathcal{Y}]\in\Gamma E$,
$\forall\mathcal{X},\mathcal{Y}\in\Gamma E$. Then, necessarily, one
also has $[\mathcal{X},\mathcal{Z}]\in\Gamma E^\perp$,
$\forall\mathcal{Z}\in\Gamma E^\perp$ \cite{V-iso}. A
weak-Hamiltonian dynamical system is a vector field $X_H$ such that
$(X_H,dH)\in E^\perp$ $(H\in C^\infty(M))$. We obtain the explicit
expression of $X_H$ and of the integrability conditions of $E$ under
the regularity condition $dim(pr_{T^*M}E)=const.$ We show that the
port-controlled, Hamiltonian systems (in particular, constrained
mechanics) \cite{{BR},{DS}} may be interpreted as weak-Hamiltonian
systems. Finally, we give reduction theorems for weak-Hamiltonian
systems and a corresponding corollary for constrained mechanical
systems.}
\end{minipage} \end{center}
\vspace*{5mm}
\section{Big-isotropic structures}
In this section we recall some basic facts concerning the
big-isotropic structures that were studied in our paper
\cite{V-iso}. All the manifolds and mappings are of class $C^\infty$
and we use the standard notation of Differential Geometry, e.g.,
\cite{KN}. In particular, $M$ is an $m$-dimensional manifold,
$\chi^k(M)$ is the space of $k$-vector fields, $\Omega^k(M)$ is the
space of differential $k$-forms, $\Gamma$ indicates the space of
global cross sections of a vector bundle, $X,Y,..$ are either
contravariant vectors or vector fields, $\alpha,\beta,...$ are
either covariant vectors or $1$-forms, $d$ is the exterior
differential and $L$ is the Lie derivative.

The vector bundle $T^{big}M=TM\oplus T^*M$ is called the {\it big
tangent bundle}. It has the natural, non degenerate metric of zero
signature (neutral metric)
\begin{equation}\label{gFinC}
g((X,\alpha),(Y,\beta))=\frac{1}{2}(\alpha(Y)+\beta(X)),
\end{equation} the non degenerate, skew-symmetric $2$-form
\begin{equation}\label{omegainC}
\omega((X,\alpha),(Y,\beta))= \frac{1}{2}(\alpha(Y)-\beta(X))
\end{equation} and the Courant bracket of cross sections
\cite{C}
\begin{equation}\label{crosetC} [(X,\alpha),(Y,\beta)] = ([X,Y],
L_X\beta-L_Y\alpha+\frac{1}{2}d(\alpha(Y)-\beta(X))).\end{equation}
\begin{defin}\label{defbigiso} {\rm A $g$-isotropic subbundle
$E\subseteq T^{big}M$ of rank $k$ $(0\leq k\leq m)$ is called a {\it
big-isotropic structure} on $M$. A big-isotropic structure $E$ is
{\it integrable} if $\Gamma E$ is closed by the Courant bracket
operation.}\end{defin}

From the properties of the Courant bracket (axiom (v) of the
definition of a Courant algebroid \cite{LWX}, see \cite{V-iso}) it
follows that if
$$(X,\alpha)\in\Gamma(E),\,(Y,\beta)\in\Gamma(E),\,
(Z,\gamma)\in\Gamma(E'),$$ where $E'=E^\perp$ is the
$g$-orthogonal bundle of $E$, then
$$g([(X,\alpha),(Z,\gamma)],(Y,\beta)) +
g((Z,\gamma),[(X,\alpha),(Y,\beta)])=0,$$ whence we see that the
integrability of $E$ is equivalent with the property that
$[E,E']\subseteq E'$ (Courant bracket).

The big-isotropic structures are a generalization of the (almost)
Dirac structures which are obtained if $k=m$. The reader can find
many examples in \cite{V-iso}, in particular the following one which
we will use later.
\begin{example}\label{exgraphP} {\rm Let $\Sigma$ be a subbubdle of rank $k$
of $T^*M$ and $P\in\chi^2(M)$ a bivector field. Then
\begin{equation}\label{eqEP} E_P=graph(\sharp_P|_{\Sigma}) =
\{(\sharp_P\sigma=i(\sigma)P,\sigma)\,/\,\sigma\in
\Sigma\}\end{equation} is a big-isotropic structure on $M$ with the
$g$-orthogonal bundle
\begin{equation}\label{eqE'P} E'_P=\{(\sharp_P\beta+Y,\beta)\,/\,\beta\in T^*M,
Y\in S=ann\,\Sigma\}.\end{equation} The structure (\ref{eqEP}) is
integrable iff \cite{V-iso}: 1) $\Sigma$ is closed with respect to
the bracket of $1$-forms defined by
\begin{equation}\label{croset1forme} \{\alpha,\beta\}_P=
L_{\sharp_P\alpha}\beta-L_{\sharp_P\beta}\alpha
-d(P(\alpha,\beta)),\end{equation} 2) the Schouten-Nijenhuis bracket
$[P,P]$ (e.g., \cite{V-carte}) satisfies the condition
\begin{equation}\label{cond2-a} [P,P](\sigma_1,\sigma_2,\beta)=0,
\hspace{3mm}\forall\sigma_1,\sigma_2\in \Sigma,\forall\beta\in
T^*M.\end{equation}}\end{example}

We may also define a big-isotropic structure on a vector space (or a
vector bundle) $V$ as an isotropic subspace $E\subseteq V\oplus
V^*$. Then, we get the subspaces $U_E=pr_VE,U_{E'}=pr_VE'$ and a
bilinear mapping $\varpi:U_E\times U_{E'}\rightarrow\mathds{R}$
given by
\begin{equation}\label{varpi} \varpi(v_1,v_2)=
\omega((v_1,a_1),(v_2,a_2))= a_1(v_2)=-a_2(v_1),\end{equation} where
$(v_1,a_1)\in E,(v_2,a_2)\in E'$ (the equalities hold and the
result is independent of the choice of $a_1,a_2$ because
$(v_1,a_1)\perp_g(v_2,a_2)$). The following result is Proposition
2.1 plus formula (2.17) of \cite{V-iso}:
\begin{prop}\label{elemechiv} For any pair of subspaces
$U_E\subseteq U_{E'}\subseteq V$ and any bilinear mapping
$\varpi:U_E\times U_{E'} \rightarrow\mathds{R}$ with a
skew-symmetric restriction to $U_E\times U_E$, there exists a
unique, big-isotropic subspace $E\subseteq V\oplus V^*$ such that
$U_E=pr_VE,U_{E'}=pr_VE'$ and $\varpi$ is the mapping
{\rm(\ref{varpi})}. The space $E$ and the orthogonal space $E'$ are
given by
\begin{equation}\label{eqE} \begin{array}{l}
E=\{(v,a)\,/\, v\in U_E,\,\forall w\in U_{E'},\,a(w)=\varpi(v,w)\},
\vspace{2mm}\\ E'=\{(w,b)\,/\, w\in U_{E'}, \,\forall v\in
U_E,\,b(v)=-\varpi(v,w)\}.
\end{array}\end{equation} The dimensions of the spaces above satisfy
the following equalities
\begin{equation}\label{reldim} dim\,E=dim\,U_E+dim \,ann\,U_{E'},\;
dim\,E'=dim\,U_{E'}+dim \,ann\,U_{E}.\end{equation}\end{prop}
\section{Weak-Hamiltonian vector fields}
The aim of this paper is to show that the big-isotropic structures
are interesting for physics and control theory because they define a
Hamiltonian formalism that may be used in applications.
\begin{defin}\label{defiweakH} {\rm \cite{V-iso}
A function $f\in C^\infty(M)$ is a {\it Hamiltonian}, respectively
{\it weak-Hamiltonian}, function if there exists a vector field
$X_f\in\chi^1(M)$ such that $(X_f,df)\in\Gamma E$, respectively
$(X_f,df)\in\Gamma E'$. The vector field $X_f$ is a {\it
Hamiltonian}, respectively {\it weak-Hamiltonian}, vector field of
$f$.}\end{defin}

The vector field $X_f$ is not unique; in the Hamiltonian case $X_f$
is defined up to the addition of any $Z\in ann\,pr_{T^*M}E'$ and in
the weak-Hamiltonian case up to $Z\in ann\,pr_{T^*M}E$. We denote by
$C^\infty_{ham}(M,E)$ the set of Hamiltonian functions, by
$C^\infty_{wham}(M,E)$ the set of weak-Hamiltonian functions and by
$\chi_{ham}(M,E),\chi_{wham}(M,E)$, respectively, the sets of
Hamiltonian and weak-Hamiltonian vector fields. It follows that
$Z\in\chi_{ham}(M,E)$ is Hamiltonian, for two functions $f_1,f_2$
iff $df_2-df_1\in ann\,U_{E'}$ and $Z\in\chi_{wham}(M,E)$ is
weak-Hamiltonian for $f_1,f_2$ iff $df_2-df_1\in ann\,U_E$.

Furthermore, if $f\in C^\infty_{ham}(M,E)$ and $h\in
C^\infty_{wham}(M,E)$ the following bracket is well defined
\begin{equation}\label{crosetfh}
\{f,h\}=-\varpi(X_f,X_h)=X_{f}h=-X_{h}f\end{equation}
and does not depend on the choice of the Hamiltonian vector fields
of the function $f,h$. The bracket (\ref{crosetfh}) is called the
{\it Poisson bracket} of the two functions.

Even though it is defined in the general case, the Poisson bracket
has interesting properties if $E$ is an integrable, big-isotropic
structure, which we assume for the moment. Then, formula
(\ref{crosetC}) shows that $\{f,h\}\in C^\infty_{wham}(M,E)$ and one
of its weak-Hamiltonian vector fields is $[X_f,X_h]$. If both
$f,h\in C^\infty_{ham}(M,E)$, their Poisson bracket is skew
symmetric and belongs to $C^\infty_{ham}(M,\\E)$. Furthermore, the
Poisson bracket satisfies the {\it Leibniz rule}
\begin{equation}\label{Leibniz} \{l,\{f,h\}\} =
\{\{l,f\},h\}\}+\{f,\{l,h\}\}, \end{equation} $\forall l,f\in
C^\infty_{ham}(M,E),h\in C^\infty_{wham}(M,E)$. Property
(\ref{Leibniz}) restricts to the Jacobi identity on
$C^\infty_{ham}(M,E)$. Thus, $C^\infty_{ham}(M,E)$ with the Poisson
bracket is a Lie algebra and $C^\infty_{wham}(M,E)$ is a module over
this Lie algebra. Also, $\chi_{ham}(M,E)$ is a Lie subalgebra of
$\chi^1(M)$ and $\chi_{wham}(M,E)$ is a module over the former for
the usual Lie bracket of vector fields.

In what follows integrability will hold only if explicitly
postulated. In the remaining part of this section we discuss some
big-isotropic structures where one has an explicit expression of a
weak-Hamiltonian vector field, a fact that is important in
applications. For instance, for a big-isotropic structure of the
form (\ref{eqEP}) formula (\ref{eqE'P}) provides such an expression:
\begin{equation}\label{wHinexP} X_H=\sharp_PdH+Y,\;
\forall Y\in S,\,H\in C^\infty(M).\end{equation} The following
proposition extends a result given in \cite{DS} for almost Dirac
structures.
\begin{prop}\label{propcazreg} Let $E$ be a big-isotropic
structure on $M$ such that
\begin{equation}\label{*reg} dim(pr_{T_x^*M}E_x)=const.\;(x\in M).
\end{equation}
Then, there exist bivector fields $\Pi\in\chi^2(M)$ such that if $H$
is a Hamiltonian, respectively a weak-Hamiltonian, function the
formulas
\begin{equation}\label{vectHinreg} X_H=\sharp_\Pi dH +
Z,\;\;Z\in ann\,pr_{T_x^*M}E',\end{equation} respectively
\begin{equation}\label{vectwHinreg}
X_H=\sharp_\Pi dH + W,\;\;W\in ann\,pr_{T_x^*M}E,\end{equation}
define the Hamiltonian, respectively weak-Hamiltonian, vector fields
of $H$.\end{prop}
\begin{proof} For a simpler notation put
\begin{equation}\label{SS'} \begin{array}{l}
\Sigma=pr_{T^*M}E,\,S=ann\,\Sigma=(TM\oplus0)\cap
E',\vspace{2mm}\\\Sigma'=pr_{T^*M}E',
\,S'=ann\,\Sigma'=(TM\oplus0)\cap E;\end{array}\end{equation} notice
that
\begin{equation}\label{SrelS'} \Sigma\subseteq \Sigma',\,S'\subseteq
S.\end{equation} We shall use Proposition \ref{elemechiv} for the
fibers $E_x$ of $E$ $(x\in M)$ taking $V=T^*_xM$ and denoting the
corresponding bilinear mapping $\varpi$ by $P_x:\Sigma_x\times
\Sigma'_x\rightarrow\mathds{R}$. Then, after changing the order of
the terms of a pair, formulas (\ref{eqE}) become
\begin{equation}\label{eqEcuP} \begin{array}{l}
E_x=\{(X,\alpha)\,/\,\alpha\in
\Sigma_x,\,\beta(X)=P_x(\alpha,\beta),\,
\forall\beta\in \Sigma'_x\},\vspace{2mm}\\
E'_x=\{(Y,\beta)\,/\,\beta\in
\Sigma'_x,\,\alpha(Y)=-P_x(\alpha,\beta),\, \forall\alpha\in
\Sigma_x\},\end{array}\end{equation} where
\begin{equation}\label{Px}
P_x(\alpha,\beta)=\frac{1}{2}(\alpha(Y)-\beta(X)),\end{equation} for
any choice of $X,Y$ such that $(X,\alpha)\in E_x,(Y,\beta)\in E'_x$
and the result is independent of this choice. Hypothesis
(\ref{*reg}), which will be called the $*$-{\it regularity
condition}, together with formulas (\ref{reldim}), show that
$\Sigma,\Sigma'$ are subbundles of $T^*M$. Therefore, we may choose
bundle decompositions
\begin{equation}\label{descaux} \Sigma'=\Sigma\oplus Q_1,\;
T^*M=\Sigma\oplus Q_1\oplus Q_2. \end{equation} Then, we can extend
$P$ to a bivector field $\Pi$ by means of the formula
\begin{equation}\label{Pi}
\Pi(\lambda,\mu)= P(\lambda',\mu')+P(\lambda',\mu'')
-P(\mu',\lambda'')\;\;(\lambda,\mu\in T^*M),\end{equation} where $'$
and $''$ denote the first and second projection in the decomposition
(\ref{descaux}) of $T^*M$, and the expressions (\ref{eqEcuP}) become
\begin{equation}\label{reconstr2} \begin{array}{l}
E_x=\{(X,\alpha)\,/\,\alpha\in \Sigma_x,\,X|_{\Sigma'_x}=
(\sharp_{\Pi_x}\alpha)|_{\Sigma'_x}\},\vspace{2mm}\\
E'_x=\{(Y,\beta)\,/\,\sharp_{\Pi_x}\beta-Y\in S_x,\,\beta\in
\Sigma'_x\}.
\end{array} \end{equation} The required formulas
(\ref{vectHinreg}), (\ref{vectwHinreg}) are a straightforward
consequence of (\ref{reconstr2}). \end{proof}

It is obvious that, in fact, only the values of the mapping $P$
actually appear in the expressions of the vector fields
(\ref{vectHinreg}), (\ref{vectwHinreg}) and two bivector fields
$\Pi_1,\Pi_2$ produce the same values $X_H$ iff they have the same
restriction $P$ to $\Sigma\times \Sigma'$. Notice also that the
formulas (\ref{vectwHinreg}) and (\ref{wHinexP}) differ only by the
fact that the former includes the restriction $dH\in \Sigma'$. In
view of (\ref{SS'}), if $(TM\oplus0)\cap E=0$ this restriction is
void, therefore, any function $H\in C^\infty(M)$ is a
weak-Hamiltonian function and formulas (\ref{eqEP}), (\ref{eqE'P})
with $P$ replaced by $\Pi$ hold. Still, $\Pi$ is not uniquely
defined.
\begin{rem}\label{obshamilt} {\rm It is always possible to consider an
arbitrary ``Hamiltonian function" $H\in C^\infty(M)$, then restrict
to the subset of the points of $M$ where $dH\in \Sigma'$
\cite{BR}.}\end{rem}

The following proposition yields the integrability conditions of a
$*$-regular, big-isotropic structure.
\begin{prop}\label{integrabinreg} Let $E$ be a $*$-regular,
big-isotropic structure with the associated subbundles
$\Sigma,S,\Sigma',S'$ and let $\Pi\in\chi^2(M)$ be such that
formulas  {\rm(\ref{vectHinreg})}, {\rm(\ref{vectwHinreg})} hold.
Then
\begin{equation}\label{formaEreg} \begin{array}{l}
E=\{(\sharp_\Pi\alpha+Z,\alpha)\,/\,\alpha\in \Sigma,Z\in S'\},
\vspace{2mm}\\
E'=\{(\sharp_\Pi\beta+W,\beta)\,/\,\beta\in \Sigma',W\in S\}.
\end{array}\end{equation}
The structure $E$ is integrable iff the following conditions are
satisfied: \\

1) the distribution $S'$ is integrable and $S$ is projectable to
the space of leaves of $S'$ (see Section 5 of {\rm\cite{V-iso}}
for this notion of projectability);\\

2) the subbundle $\Sigma$ is closed by the $\Pi$-brackets
{\rm(\ref{croset1forme})} and $\forall\alpha\in\Gamma\Sigma$,
$\forall\beta\in\Gamma\Sigma'$ one has $\{\alpha,\beta\}\in\Gamma\Sigma'$;\\

3) $\forall\alpha_1,\alpha_2\in\Gamma \Sigma,\beta\in\Gamma \Sigma'$
one has $$[\Pi,\Pi](\alpha_1,\alpha_2,\beta)=0.$$
\end{prop}
\begin{proof}
The formulas (\ref{formaEreg}) have the same proof like
(\ref{vectHinreg}), (\ref{vectwHinreg}).

If we use the Gelfand-Dorfman formula
\begin{equation}\label{GD}
\Pi(\{\alpha_1,\alpha_2\}_{\Pi},\beta) =
\gamma([\sharp_\Pi\alpha_1,\sharp_\Pi\alpha_2]) +\frac{1}{2} [\Pi,
\Pi](\alpha_1,\alpha_2,\beta),
\end{equation} we get
\begin{equation}\label{crosetreg} \begin{array}{l}
[(\sharp_\Pi\alpha_1+Z_1,\alpha_1),(\sharp_\Pi\alpha_2+Z_2,\alpha_2)]\vspace{2mm}\\
=(\sharp_\Pi(\{\alpha_1,\alpha_2\}_\Pi-L_{Z_2}\alpha_1+L_{Z_1}\alpha_2)
-\sharp_{L_{Z_2}\Pi}\alpha_1+\sharp_{L_{Z_1}\Pi}\alpha_2\vspace{2mm}\\
+[Z_1,Z_2]-\frac{1}{2}i(\alpha_1\wedge\alpha_2)[\Pi,\Pi],\vspace{2mm}
\{\alpha_1,\alpha_2\}_\Pi-L_{Z_2}\alpha_1+L_{Z_1}\alpha_2),\end{array}\end{equation}
where $Z_1,Z_2\in S',\alpha_1,\alpha_2\in\Sigma$.

The structure $E$ is integrable iff the right hand side of
(\ref{crosetreg}) belongs to $E$ and we may brake this condition
into the cases: a) $\alpha_1=\alpha_2=0$, b) $Z_1=Z_2=0$, c)
$Z_1=0,\alpha_2=0$ (equivalently $Z_2=0,\alpha_1=0$).

In case a) the condition becomes $([Z_1,Z_2],0)\in E$, which is
equivalent to the first assertion of condition 1) of the
proposition.

Furthermore in case b) the bracket (\ref{crosetreg}) belongs to
$E$  iff the first assertion of condition 2) and condition 3)
hold.

Finally, a technical computation shows that if $\alpha\in\Sigma,
\beta\in\Sigma', Z\in S'\subseteq S$ then
\begin{equation}\label{eqaux1}
L_{Z}\Pi(\alpha,\beta)=\{\alpha,\beta\}_\Pi(Z).\end{equation}

Now, in case c) the right hand side of (\ref{crosetreg}) is
\begin{equation}\label{auxnou}
-(\sharp_\Pi(L_{Z_2}\alpha_1)+\sharp_{L_{Z_2}\Pi}\alpha_1,L_{Z_2}\alpha_1),
\end{equation}
which belongs to $E$ iff
$$L_{Z_2}\alpha_1\in\Sigma,\,\sharp_{L_{Z_2}\Pi}\alpha_1\in S'.$$
From (\ref{SS'}) and (\ref{eqaux1}) it follows that the two
conditions mentioned above are equivalent with the second assertions
of 1) and 2), respectively. \end{proof}
\begin{rem}\label{obsproj}
{\rm Let $E$ be an integrable, $*$-regular, big-isotropic structure.
Then, Corollary 5.1 of \cite{V-iso} shows that $E$ is projectable
with respect to the foliation $S'$, and the projection of $E$ onto
the local spaces of the slices of $S'$ is an integrable,
big-isotropic structure of the type discussed in Example
\ref{exgraphP}.}\end{rem}
\section{Port-controlled dynamical systems}
In this section we present some applications where
weak-Hamiltonian vector fields can be used. Following \cite{DS}, a
physical network is a sum of {\it port-controlled, generalized,
Hamiltonian systems} with interconnections. Many concrete
examples, in particular constrained mechanics, are discussed in
\cite{{BR},{DS}}. We shall give weak-Hamiltonian interpretations
of such port-controlled systems.

With the notation of \cite{DS}, a port-controlled, generalized,
Hamiltonian system is a system of equations of the following form
\begin{equation}\label{Schaft}
\begin{array}{l}
\dot{x}=J(x)\frac{\partial H}{\partial x}(x)+g(x)f,\vspace{2mm}\\
e=g^T(x)\frac{\partial H}{\partial x}(x),\end{array}\end{equation}
where a dot denotes time-derivative and one uses the matrix
notation. In (\ref{Schaft}) $x=(x_i)$ $(i=1,...,n)$ is the column of
energy variables, which are local coordinates on a manifold $N$ seen
as the phase space, $H$ is the total stored energy, $J$ is a
skew-symmetric $(n,n)$-matrix, $f=(f_j)\in\mathds{R}^p$
$(j=1,...,p)$ is the column of flows, $g$ is an $(n,p)$-matrix,
$e=(e_j)$ is the column of efforts and $T$ denotes matrix
transposition.

The evolution of the system is defined by the differential equations
on the first line of (\ref{Schaft}) where a choice of functions
$f_j=f_j(x)$ is made. If we see $J$ as a bivector field on $N$ and
$g$ as a vector bundle morphism $g:N\times\mathds{R}^p\rightarrow
TN$, these differential equations are equivalent with the
weak-Hamiltonian vector field
\begin{equation}\label{SchaftHam} X_H=\sharp_JdH+gf\end{equation} of
the function $H$ with respect to the big-isotropic structure
\begin{equation}\label{exS1} E_J=graph(\sharp_J|_{\Sigma}),\;\;
\Sigma=ann\,S,
\end{equation} where $S$ is any distribution on $N$ such that
$im\,g\subseteq S$. If $rank\,g=const.$ and $S=im\,g$, formula
(\ref{SchaftHam}) is that of all the weak-Hamiltonian vector fields
of $H$. Since (\ref{exS1}) is of the type (\ref{eqEP}) we see that a
port-controlled system has a weak-Hamiltonian interpretation with
respect to an integrable big-isotropic structure iff there exists a
subbundle $S\subseteq TM$ that contains $im\,g$, $ann\,S$ is closed
by the bracket (\ref{croset1forme}) for $P=J$ and $J$ satisfies the
condition (\ref{cond2-a}).

Moreover, we can show that the whole system (\ref{Schaft}) may be
seen as a weak-Hamiltonian vector field on $M=N\times\mathds{R}^p$.
For this purpose, notice that $g$ defines a bivector field
$G\in\chi^2(M)$ given by
\begin{equation}\label{GSchaft} G_{(x,f)}(\alpha'+\alpha'',\beta'+\beta'')
=\beta''(g^T(x)\alpha')-\alpha''(g^T(x)\beta'),\end{equation} where
$x\in N, f\in \mathds{R}^p,\alpha',\beta'\in T^*_xN,\alpha'',\beta''
\in T^*_f\mathds{R}^p\approx\mathds{R}^p$ and
$g^T(x):T^*_xN\rightarrow\mathds{R}^{p*}$. Then, we have the
bivector field $P=J+G\in\chi^2(M)$ and the weak-Hamiltonian vector
field
\begin{equation}\label{vectorinS2} X_H=\sharp_{J+G}dH+gf
\end{equation} of $H$ with respect to any big-isotropic structure
$graph(\sharp_{J+G}|_{\Sigma})$, where $S=ann\,\Sigma$ is a regular
distribution on $N$ that contains $im\,g$. The integral lines of the
vector field (\ref{vectorinS2}) are given by (\ref{Schaft}) where
$e_j$ are the time derivatives of the coordinates of the factor
$\mathds{R}^p$ of $M$ and one uses the natural identification of
$T\mathds{R}^p$ with $\mathds{R}^p$. The integrability conditions of
$graph(\sharp_{J+G}|_{\Sigma})$ are provided by (\ref{croset1forme})
and (\ref{cond2-a}) again.

In \cite{DS} one also defines {\it port-controlled Hamiltonian
systems with constraints}, which have the form
\begin{equation}\label{Schaftc}
\begin{array}{l}
\dot{x}=J(x)\frac{\partial H}{\partial x}(x)+g(x)f+b(x)\lambda,\vspace{2mm}\\
e=g^T(x)\frac{\partial H}{\partial x}(x),\,0=b^T(x)\frac{\partial
H}{\partial x}(x)\end{array}\end{equation} where the notation is
like in (\ref{Schaft}), $b$ is an $(n,k)$-matrix and
$\lambda\in\mathds{R}^k$. As in the non-constrained case, the system
(\ref{Schaftc}) is a weak-Hamiltonian system on
$M=N\times\mathds{R}^{p+k}$, where the Hamiltonian function $H$ is
required to satisfy the constraint $b^T(dH)=0$.

Consider the port-controlled system (\ref{Schaft}) again. It is
called {\it energy-preserving} \cite{DS} if the vectors
$f\in\mathds{R}^p,e\in\mathds{R}^p\approx(\mathds{R}^p)^*$ are
assumed to satisfy the condition $(f,e)\in\Delta(x)$ where
$\Delta(x)$ is a maximal (i.e., $p$-dimensional), isotropic subspace
of $\mathds{R}^p\times\mathds{R}^{*p}$ parameterized by $x\in N$.
The reason for this name is that, then, the energy $H$ is preserved
along the integral lines of the vector field (\ref{SchaftHam}) of
the system. Indeed, in view of the second equation (\ref{Schaft})
and since $(f,e)\in\Delta$ implies $e(f)=0$, we have
$$\dot{H}=X_HH=0.$$

Then, it turns out that the differential equations of the first line
of (\ref{Schaft}) are equivalent with a Hamiltonian vector field
with respect to an almost Dirac structure. We give a more conceptual
proof of this result proven differently in Proposition 2.2 of
\cite{DS}.

With the notation of (\ref{SchaftHam}), put
\begin{equation}\label{2.2DS}
D=\{(\sharp_J\alpha+gf,\alpha)\,/\,(f,g^T\alpha)\in\Delta\}
\subseteq TN\oplus T^*N.\end{equation} The isotropy of $\Delta$
implies that $D$ is a big-isotropic structure on $N$ and we shall
compute $dim\,D$ for any fixed point $x\in N$. Denote
$\Delta'=\Delta\cap(\mathds{R}^p\times im\,g^T)$. Then the
correspondence
$$(\sharp_J\alpha+gf,\alpha)\mapsto(f,g^T\alpha)$$ produces a
surjective homomorphism
\begin{equation}\label{phi} \phi:D\rightarrow \Delta'/
\Delta'\cap(ker\,g\oplus0)\end{equation} with
$$ker\,\phi=\{(\sharp_J\alpha,\alpha)\,/\,g^T\alpha=0\},$$ whence,
\begin{equation}\label{dimnucleu}dim\,ker\,\phi=n-rank\,g\;(n=dim N).
\end{equation}

On the other hand, if we notice that
$$\mathds{R}^p\times im\,g^T=(ker\,g)^\perp$$ (perpendicularity is
with respect to the neutral metric of
$\mathds{R}^p\oplus\mathds{R}^{p*}$ and the result holds because
the two spaces are orthogonal and the sum of their dimensions is
$2p$), we get
$$\Delta'=\Delta\cap(ker\,g)^\perp= \Delta^\perp\cap(ker\,g)^\perp
=(\Delta+ker\,g)^\perp$$ ($\Delta^\perp=\Delta$ because of the
maximal isotropy of $\Delta$). Now, if $dim(\Delta\cap ker\,g)=i$
the known formula
$$dim(\Delta+ker\,g)=dim\Delta+dim\,ker\,g-i,$$ implies
$$dim\Delta'=2p-dim(\Delta+ker\,g)=rang\,g+i.$$

Together with (\ref{phi}) and (\ref{dimnucleu}), the previous result
gives $dim\,D=(n-rank\,g)+[(rang\,g+i)-i]=n$, hence, $D$ is an
almost Dirac structure. Furthermore, for $X_H$ given by
(\ref{SchaftHam}) and since we asked that $(f,e)\in\Delta$, we have
$(X_H,dH)\in\Gamma\,D$ and $X_H$ is a Hamiltonian vector field of
$H$.
\begin{rem}\label{obsptnetw}{\rm
The systems discussed in \cite{{BR},{DS}} are direct sums of
port-controlled systems on a product manifold where the components
may not be energy preserving but the sum  is such. These are {\it
energy-preserving physical networks} and the corresponding $\Delta$
is a {\it power-preserving interconnection} between the
port-controlled components \cite{DS}.}\end{rem}
\begin{rem}\label{nedif} {\rm The structure (\ref{2.2DS}) may
present a technical difficulty: even if $\Delta(x)$ is
differentiable with respect to $x\in N$, $D$ may not be
differentiable. For instance, if $\Delta(x)=\mathds{R}^p\oplus0$
one has
\begin{equation}\label{DptexDS}
D=\{(\sharp_J\alpha+Z,\alpha)\,/\,\alpha\in ann\,im\,g,\, Z\in
im\,g\}\end{equation} and $D$ is not differentiable if $rang\,g$ is
not constant. If $rang\,g=const.$, (\ref{DptexDS}) has the same form
as $E$ of (\ref{formaEreg}), with $E'=E=D$, and the integrability
conditions will be like in Proposition \ref{integrabinreg}, i.e., :
1) $im\,g$ is integrable, 2) $ann\,im\,g$ is closed by the
$J$-bracket of $1$-forms, 3) $[J,J]|_{ann\,im\,g}=0$.}\end{rem}

It was shown in \cite{{BR},{DS}} that the dynamical systems
(\ref{Schaft}) include the constrained mechanical systems. Here, we
give a straightforward, weak-Hamiltonian interpretation of a
constrained mechanical system.

A mechanical system has a configuration space, which is a manifold
$Q$, the space of the velocities, which is the tangent bundle $TQ$,
and the space of the phases, which is the cotangent bundle $T^*Q$.
Constraints consist of a $k$-dimensional distribution $L$ on $Q$. In
Hamiltonian mechanics, the differential equations of the motion are
those of the integral lines of a vector field of the form
\begin{equation}\label{Hlegaturi}X=\sharp_P
dH+\sharp_P(\pi^*\alpha)\in\chi^1(T^*Q),\end{equation} where $P$ is
defined by $\sharp_P\circ\flat_\omega=-Id$, $\omega$ being the
canonical symplectic form of $T^*Q$ (\cite{MR}, Section 6.2), $H$ is
the Hamiltonian of the system, $\pi:T^*Q\rightarrow Q$ is the
natural projection and $\alpha\in ann\,L$ (e.g., \cite{BR}).

The constraint distribution $L$ produces a natural,
$\omega$-isotropic subbundle
\begin{equation}\label{S-L} S_L=\{\sharp_P(\pi^*\alpha)
\,/\,\alpha\in ann\,L\}\subseteq T(T^*Q).\end{equation}
The corresponding $\omega$-orthogonal subbundle is
\begin{equation}\label{SperpL}
S_L^{\perp_\omega}=ann(\pi^*(ann\,L))=\{\mathcal{X}\in T(T^*Q)\,/\,
\pi_*\mathcal{X}\in L\}.\end{equation}

A comparison with formula (\ref{wHinexP}) shows that the vector
field (\ref{Hlegaturi}) is weak-Hamiltonian with respect to the big
isotropic structure $E_L=graph(\sharp_P|_{ann\,S_L})$.

The structure $E_L$ is integrable iff $ann\,S_L$ is closed by the
bracket (\ref{croset1forme}); the fact that $\omega$ is a symplectic
form implies the Poisson condition $[P,P]=0$, hence, (\ref{cond2-a})
holds too. Notice that $\sigma\in ann\,S_L$ is equivalent with
$\sharp_P\sigma\in S_L^{\perp_\omega}$ and, since \cite{V-carte}
$$\sharp_P\{\sigma_1,\sigma_2\}_P=[\sharp_P\sigma_1,\sharp_P\sigma_2],$$
it follows that $E_L$ is integrable iff the distribution
$S^{\perp_\omega}_L$ is integrable. Now, let us recall that $L$
itself is integrable iff, $\forall\alpha\in ann\,L$, $d\alpha$
belongs to the ideal spanned by $ann\,L$. Since
$ann\,S^{\perp_\omega}_L=\pi^*(ann\,L)$ and $\pi^*$ is injective,
the same condition characterizes the integrability of $S_L^\perp$.
Therefore, like in the Dirac interpretation of
\cite{BR}, the structure $E_L$ is integrable iff $L$ is integrable,
i.e., iff the system has {\it holonomic constraints}.
\section{Symmetries  and reduction} In this section we
extend some results on symmetries and reduction from Hamiltonian to
weak-Hamiltonian systems. The case of Hamiltonian systems on a Dirac
manifold was treated in \cite{{BR},{BS}}.
\begin{defin}\label{sym} {\rm A vector field $Z\in\chi^1(M)$ is an
{\it infinitesimal symmetry} of a big-isotropic structure $E$ if
\begin{equation}\label{infsym} (L_ZX,L_Z\alpha)\in\Gamma E,\hspace{3mm}
\forall(X,\alpha)\in\Gamma E.\end{equation} A diffeomorphism
$\varphi:M\rightarrow M$ is a {\it symmetry} of $E$ if
\begin{equation}\label{symglobal} (\varphi_*X,\varphi^{*-1}\alpha)\in\Gamma
E,\hspace{3mm} \forall(X,\alpha)\in\Gamma
E.\end{equation}}\end{defin}

Obviously, the flow of an infinitesimal symmetry consists of
symmetries of $E$. Furthermore, for (infinitesimal) symmetries the
conditions required for $E$ also hold for the $g$-orthogonal space
$E'$ of $E$ because the neutral metric $g$ is invariant by any
(infinitesimal transformation) diffeomorphism of $M$.
\begin{prop}\label{symreg} Let $E$ be a $*$-regular, big-isotropic
structure defined by formulas {\rm(\ref{formaEreg})}. 1). The
diffeomorphism $\varphi:M\rightarrow M$ is a symmetry of $E$ iff the
subbundles $S,S'$ are invariant by $\varphi_*$ and for all $\beta\in
\Sigma'$ one has
$\sharp_{\varphi_*(\Pi\circ\varphi^{-1})}\beta=\sharp_\Pi\beta$. 2).
The vector field $Y\in\chi^1(M)$ is an infinitesimal symmetry of $E$
iff $\forall Z\in\Gamma S'$, $\forall\alpha\in\Gamma \Sigma$,
$\forall\beta\in\Gamma \Sigma'$ one has
\begin{equation}\label{symP} [Y,Z]\in\Gamma S',
\;L_Y\alpha\in\Gamma \Sigma,\; L_Y\Pi(\alpha,\beta)=0.
\end{equation}
The conditions stated in 1), 2) depend only on the mapping $P$
defined by {\rm(\ref{Px})}.\end{prop}
\begin{proof} The notation used here is that of formula
(\ref{formaEreg}).

1). If $\varphi$ is a symmetry then, for all $\alpha\in \Sigma,Z\in
S'$, we have
\begin{equation}\label{auxsym}
(\varphi_*(\sharp_\Pi\alpha+Z),\varphi^{*-1}\alpha)=
(\sharp_{\varphi_*\Pi}(\varphi^{*-1}\alpha)+\varphi_*Z,\varphi^{*-1}\alpha)
\end{equation}
$$\stackrel{(\ref{formaEreg})}{=}(\sharp_{\Pi\circ\varphi}(\varphi^{*-1}\alpha)+U,
\varphi^{*-1}\alpha),$$ where $U\in S'$. The same must hold for all
$\alpha\in \Sigma',Z\in S$ with $U\in S$ because $\varphi$ also
preserves the orthogonal subbundle $E'$. It is easy to derive 1)
from (\ref{auxsym}) and to see that 1) also is the sufficient
condition for (\ref{auxsym}) to hold.

2). From the first formula (\ref{formaEreg}), we see that $Y$ is an
infinitesimal symmetry iff for all $\alpha\in \Sigma$ one has
$L_Y\alpha\in \Sigma$ and $\sharp_{L_Y\Pi}\alpha+[Y,Z]\in\Gamma S'$.
By looking at the cases $\alpha=0$ and $Z=0$ separately we get the
required conclusion.

The last assertion of the proposition is obvious.
\end{proof}

An infinitesimal symmetry $Z$ acts on Poisson brackets as a
derivation. Indeed, take $f\in C^\infty_{ham}(M,E)$, $h\in
C^\infty_{wham}(M,E)$ and corresponding pairs $(X_f,df)\in\Gamma
E,(X_h,dh)\in\Gamma E'$. From (\ref{infsym}), it follows that
$$([Z,X_f],d(Zf))\in\Gamma E,\; ([Z,X_h],d(Zh))\in\Gamma E',$$
whence $$Z\{f,h\}=Z(X_fh)=[Z,X_f]h+X_f(Zh)=\{Zf,h\}+\{f,Zh\}.$$

If the structure $E$ is integrable, any Hamiltonian vector field
$Z\in\chi_{ham}(M,\\ E)$ is an infinitesimal symmetry. Indeed,
assume that $Z=X_f$, $f\in C^\infty(M)$, and $(X,\alpha)\in\Gamma
E$. The integrability of $E$ implies
$$[(X_f,df),(X,\alpha)]=([X_f,X],L_{X_f}\alpha)\in\Gamma E,$$
which is the required symmetry property.

Let $H$ be a weak-Hamiltonian function on $(M,E)$. Then, we are
interested in $H$-preserving, infinitesimal and global symmetries,
i.e., vector fields $Z$ that satisfy (\ref{infsym}) and $ZH=0$, and
diffeomorphisms $\varphi:M\rightarrow M$ that satisfy
(\ref{symglobal}) and $H\circ\varphi=H$. The following proposition
is in the spirit of Noether's theorem \cite{{BR},{MR}}.
\begin{prop} Let $E$ be an integrable, big-Hamiltonian structure on $M$.
A Hamiltonian vector field $Z\in\chi_{ham}(M,E)$ is an
$H$-preserving infinitesimal symmetry for $H\in
C^\infty_{wham}(M,E)$ iff $Z$ is the Hamiltonian vector field of a
function $f$ such that $\{f,H\}=0$.\end{prop}
\begin{proof} We already know that $Z$ is an infinitesimal symmetry.
Then, the orthogonality of the pairs $(Z,df),(X_H,dH)$ gives
$$ZH-\{f,H\}=0$$ and this shows the equivalence between
$ZH=0$ and the condition required by the proposition.
\end{proof}

We may define a {\it first integral} of a weak-Hamiltonian
dynamical system $X_H$ to be a function $f\in C^\infty_{ham}(M,E)$
such that $\{f,H\}=0$ But, the usual properties of first integrals
hold only in the integrable case; then, the Hamiltonian vector
fields $X_f$ of the first integral $f$ are $H$-preserving
infinitesimal symmetries and the Poisson bracket of two first
integrals of $X_H$ is a first integral again because of the
Leibniz property (\ref{Leibniz}).

Now, let us refer to reduction. In \cite{V-iso} we discussed the
reduction of a big-isotropic structure $E$ on $M$ and we recall
the main results. Let $\iota:N\rightarrow M$ be an embedded
submanifold of $M$. Then, the formula
\begin{equation}\label{pullback} \iota^*(E_x) = \{(X,\iota^*\alpha)\,/\,
X\in T_xN,\,\alpha\in T_x^*M,\,(X,\alpha)\in E_x\},\end{equation}
where $x\in N$, defines the pullback $\iota^*E$ of $E$ to $N$.
$\iota^*E$ is a field of big-isotropic subspaces of $T^{big}(M)$
and, if this field is a differentiable subbundle of $T^{big}(M)$,
we say that the submanifold $N$ is {\it $E$-proper} with the
induced big-isotropic structure $\iota^*E$. Moreover, if $E$ is
integrable the same holds for $\iota^*E$.
\begin{rem}\label{generalconstr} {\rm A $E$-proper submanifold
$\iota:N\rightarrow M$ of $(M,E)$ may be seen as a {\it general
constraint} and a {\it constrained weak-Hamiltonian system} may be
defined as a weak-Hamiltonian vector field $X_H$ $(H\in
C^\infty(M))$ on the manifold $(N,\iota^*E)$.}\end{rem}

Furthermore, assume that the $E$-proper submanifold $N$ of $M$ has a
foliation $\mathcal{F}$ with a paracompact, Hausdorff quotient
manifold $Q=N/\mathcal{F}$ and the natural projection
$\pi:N\rightarrow Q$. Then, the formula
\begin{equation}\label{defpush} \pi_*(\iota^*E_x)=
\{(\pi_*X,\alpha)\,/\,X\in T_xN,\,\alpha\in T^*_{\pi(x)}Q,\,
(X,\pi^*\alpha)\in \iota^*E_x\}
\end{equation} defines a big-isotropic subspace of
$T^{big}_{\pi(x)}Q$, $\forall x\in N$.

Assume that the following two {\it reducibility conditions} are
satisfied:\\ \indent R1) $T\mathcal{F}\oplus0\subseteq\iota^* E$,\\
\indent R2) every vector field $Y\in\chi^1(N)$ that is tangent to
$\mathcal{F}$ is an infinitesimal symmetry of $\iota^*E$.

Then, $E^{red}=\pi_*(\iota^*E)$ given by (\ref{defpush}) is a well
defined, big-isotropic structure on $Q$ called the {\it reduced
structure} of $E$ via $(N,\mathcal{F})$. Moreover, if $E$ is
integrable condition R1) implies R2) and the reduced structure
$E^{red}$ is integrable too \cite{V-iso}.
\begin{theorem}\label{thred} Let $E$ be a
big-isotropic structure on the manifold $M$. Assume that the
connected, Lie group $G$ acts on $M$ by symmetries of $E$ that keep
fixed an embedded submanifold $\iota:N\rightarrow M$. Assume that
the restriction of the action of $G$ to $N$ is proper and free and
denote by $\mathcal{F}$ the foliation of $N$ by the orbits of $G$.
Finally, assume that the following reducibility condition holds
\vspace{2mm}\\ \hspace*{1cm}R)\hspace{2mm} for any infinitesimal
transformation $Z$ of $G$, $\exists\alpha\in ann\,TN$\\
\hspace*{1cm}such that $(Z|_N,\alpha)\in E|_N$.\vspace{2mm}\\ Then,
there exists a Hausdorff manifold $Q=N/\mathcal{F}$ endowed with a
reduced, big-isotropic structure $E^{red}$ and if $E$ is integrable
$E^{red}$ is integrable too.\end{theorem}
\begin{proof} If $E$ is integrable, this is Corollary 5.2 of
\cite{V-iso}. But, the fact that condition R) is equivalent with R1)
holds in the non-integrable case too. Condition R2) holds for the
infinitesimal transformations $Z$ of $G$ on $N$ because of the
invariance of $E$ and $N$. This implies the fact that any vector
field spanned by such infinitesimal transformations is also an
infinitesimal symmetry of $\iota^*E$. Indeed, for any $f\in
C^\infty(N)$ and $(X,\iota^*\alpha)\in\Gamma\iota^*E$ one has
$$ (L_{fZ}X,L_{fZ}(\iota^*\alpha))=f(L_ZX,L_Z(\iota^*\alpha))-(Xf)(Z,0)
+(\iota^*\alpha)(Z)(0,df),$$ where
$(L_ZX,L_Z(\iota^*\alpha))\in\iota^*E$ because $Z$ is an
infinitesimal symmetry, $(Z,0)\in\iota^*E$ by R1), and
$(\iota^*\alpha)(Z)=0$ because the isotropy of $\iota^*E$ implies
$(X,\iota^*\alpha)\perp_g (Z,0)$. Hence, R2) holds as stated and we
are done.\end{proof}

Theorem \ref{thred} is straightforwardly enhanced by the following
result, which we call a theorem because of its in-principle
importance.
\begin{theorem}\label{redH} Assume that the notation and hypotheses of Theorem
\ref{thred} hold and that we have a $G$-invariant, weak-Hamiltonian
function $H\in C^\infty_{wham}(M,E)$ with a weak-Hamiltonian
vector field $X_H$ such that $X_H(x)\in T_xN$, $\forall x\in N$
and $X_H|_N$ is $\mathcal{F}$-projectable. Then, the function
$H|_N$ is the lift by $\pi$ of a function $H^{red}\in
C^\infty_{wham}(Q,E^{red})$ and $\pi_*(X_H|_N)$ is a
weak-Hamiltonian vector field $X^{red}_{H^{red}}\in\chi^1(Q)$ of
$H^{red}$.\end{theorem}
\begin{proof} Notice that $(X_H|_N,d(H\circ\iota))\in\iota^*E'$, where
the latter is defined like $\iota^*E$ and is equal to the
orthogonal space $(\iota^*E)'$ \cite{V-iso}. The existence of
$H^{red}$ and $X^{red}_{H^{red}}$ is obvious and (\ref{defpush})
shows that $(X^{red}_{H^{red}},dH^{red})\in(E^{red})'$.\end{proof}
\begin{rem}\label{obsptred} {\rm Each of
the following two conditions: i) $X_H$ is $G$-invariant, ii)
$\iota^*E'\cap(TN\oplus0)=T\mathcal{F}$ implies the
$\mathcal{F}$-projectability of $X_H|_N$. Under condition i), it
is obvious that $X_H|_N$ is $\mathcal{F}$-projectable.
Furthermore, if $Z\in\chi^1(M)$ is an infinitesimal action of $G$,
$Z$ is an infinitesimal symmetry of $E$ and $(X_H,dH)\in E'$
implies $(L_ZX_H,L_ZdH)=([Z,X_H],0)\in E'$. Since both $Z$ and
$X_H$ are tangent to $N$, we get $([Z,X_H]|_N,0)\in\iota^* E'$
and, if hypothesis ii) holds, $[Z,X_H]|_N\in T\mathcal{F}$.
Therefore, again, $X_H|_N$ is projectable to $Q$.}\end{rem}

Thus, we can simplify the integration of a weak-Hamiltonian,
dynamical system by reduction if we have a convenient group of
symmetries and a nice invariant submanifold.

Like for the usual Hamiltonian systems, the required submanifold may
come from a momentum map. We will say that an $E$-preserving action
of a connected Lie group $G$ on $(M,E)$ is a {\it Hamiltonian
action} if the infinitesimal transformations $Z$ of $G$ are
Hamiltonian vector fields, i.e., $\exists f\in C^\infty(M)$ such
that $(Z,df)\in\Gamma E$. Like in the Poisson case (e.g.,
\cite{V-carte}, Proposition 7.25), it follows that the action is
Hamiltonian iff it preserves $E$ and $\exists J\in
C^\infty(M,\mathcal{G}^*)$ such that
\begin{equation}\label{moment} (\xi_M,d(\xi\circ J))\in \Gamma E,\;\;
\forall\xi\in\mathcal{G},\end{equation} where $\mathcal{G}$ is the
Lie algebra of $G$ and $\xi_M$ is the infinitesimal action of $\xi$
on $M$. Such a function $J$ is a {\it momentum map}. Notice that if
$E$ is integrable and the action has a momentum map then the action
necessarily preserves $E$ because the Hamiltonian vector fields are
infinitesimal symmetries of $E$. Finally, a momentum map $J$ is
equivariant if $J(g(x))= coad_g(J(x))$, $\forall g\in G,x\in M$.

From Theorems \ref{thred} and \ref{redH}, and with the notation
there, we get
\begin{corol}\label{corolJ} Consider an action of $G$ on $M$ that
preserves $E$ and has an equivariant momentum map
$J:M\rightarrow\mathcal{G}^*$ such that $0$ is a regular value of
$J$. Assume that $G$ acts properly and freely on the $G$-invariant
submanifold $N=J^{-1}(0)\subseteq M$ giving rise to the quotient
manifold $Q=N/\mathcal{F}$ where the leaves of $\mathcal{F}$ are
the orbits of $G|_N$. Then $Q$ has the reduced, big-isotropic
structure $E^{red}$ of $E$, which is integrable if $E$ is
integrable. Furthermore, consider a pair $(X_H,dH)\in \Gamma E'$
where $H\in C^\infty_{wham}(M,E)$ is $G$-invariant and $X_H|_N$ is
$\mathcal{F}$-projectable. Then $(X_H,dH)|_N$ projects to a pair
$(X^{red}_{H^{red}},dH^{red})\in \Gamma (E^{red})'$ and one has a
reduced, weak-Hamiltonian system on $Q$.\end{corol}
\begin{proof} For the first assertion we just have to check
condition R). If $Z=\xi_M$ for $\xi\in\mathcal{G}$ then
$(Z,d(\xi\circ J))\in\Gamma E$ and, since $J$ is constant on $N$,
$d(\xi\circ J)\in ann\,TN$, which is the required condition. For the
second assertion we have to check that $X_H$ is tangent to $N$. This
holds because the invariance of $H$ implies $X_H(\xi\circ
J)=-X_{\xi\circ J}H=-\xi_MH=0$.
\end{proof}
\begin{rem}\label{corolintegr} {\rm
The $G$-invariance of $H$ is equivalent with $\{\xi\circ J,H\}=0$,
$\forall\xi\in\mathcal{G}$. Hence, like in symplectic mechanics,
given a system $X_H$ on $(M,E)$, we should look for symmetry groups
$G$ that lead to reduction by looking for first integrals $f_i$ of
$H$ such that $X_{f_i}$ are infinitesimal symmetries of $E$ and
$span\{X_{f_i}\}$ is a Lie algebra.}
\end{rem}
\begin{rem}\label{valnecrit} {\rm We can reformulate Corollary
\ref{corolJ} for an arbitrary non-critical value $\gamma$ of $J$ and
the level set $N=J^{-1}(\gamma)$. Indeed, if the group $G$ satisfies
the hypotheses of Corollary \ref{corolJ} and $G'$ is a connected
subgroup of $G$ with the Lie algebra
$i:\mathcal{G}'\subseteq\mathcal{G}$, it follows easily that
$J'=i^T\circ J$ is an equivariant momentum map of the action of $G'$
on $M$. In particular, if $G'=G_\gamma$ is the isotropy subgroup of
$\gamma\in\mathcal{G}^*$ with respect to the coadjoint action we
have $J^{-1}(\gamma)=J'^{-1}(0)$, and we may use Corollary
\ref{corolJ} for the connected component of the unit of $G_\gamma$
instead of $G$. The result will be a version of the
Marsden-Weinstein reduction theorem in the present
context.}\end{rem}

We finish by discussing the application of Corollary \ref{corolJ} to
the constrained mechanical system described at the end of Section 3,
with the notation used there, i.e., the configuration space is $Q$,
the constraint distribution is $L\subseteq TQ$ and the associated
big-isotropic structure is $E_L$. Assume that $G$ is a connected Lie
group acting on $Q$ such that the distribution $L$ is {\it strongly
invariant}, by which we mean the following two conditions: a)
$\forall g\in G$, $g_*(L)=L$, b) $\forall x\in M$,
$T_x(G(x))\subseteq L_x$ ($G(x)$ is the $G$-orbit of the point $x$).
Then, the derivative mappings yield a group $G^{tg^*}$that acts on
the phase space $T^*Q$ by symplectomorphisms of the canonical
symplectic form $\omega$ and preserves the big-isotropic structure
$E_L$. Furthermore, there exists a well known, equivariant, momentum
map $J:T^*Q\rightarrow\mathcal{G}^*$ for the symplectic structure of
$T^*Q$ defined by
$$<J(\alpha),\xi>=<\alpha,\xi_Q>\hspace{5mm}(\alpha\in
T^*Q,\xi\in\mathcal{G})$$ (e.g., \cite{MR}, Theorem 12.1.4).

The fact that $J$ is a momentum map for $\omega$ means that we have
\begin{equation}\label{mominex}
\xi_{T^*Q}=\sharp_P d(\xi\circ
J)\hspace{3mm}(\sharp_P\flat_\omega=-Id).\end{equation} But,
condition b) of the strong invariance of $L$ also implies
$d(\xi\circ J)\in ann(\sharp_P(ann\,L))$. Indeed, $\forall\alpha\in
ann\,L$ we have
$$<d(\xi\circ J),\sharp_P(\pi^*\alpha> \stackrel{(\ref{mominex})}{=}
-\pi^*\alpha(\xi_{T^*Q})=-\alpha(\xi_Q),$$ which vanishes because of
b). Thus, $J$ also is a momentum map with respect to the structure
$E_L$ and we get
\begin{corol}\label{redconstrmech} Let $(Q,L)$ be a constrained
mechanical system with the Hamiltonian function $H$ and the
Hamiltonian vector field $X_H$. Assume that the connected Lie
group $G$ acts on $Q$ such that: 1) $G$ strongly preserves $L$, 2)
$G^{tg*}$ preserves the pair $(H,X_H)$. Let
$J:T^*Q\rightarrow\mathcal{G}^*$ be the naturally associated
momentum map and assume that $0$ is a regular value of $J$ and
that the orbits of $G^{tg*}|_N$ are the leaves of a foliation
$\mathcal{F}$ of $N$ by the leaves of a submersion
$\pi:N\rightarrow Q$, where $Q$ is a Hausdorff, differentiable
manifold. Then the system admits a reduction to $Q$ via
$(N,\mathcal{F})$.\end{corol}

Notice that the constraints may be non-holonomic.\\

{\it Acknowledgement}. Part of the work on this paper was done
during the author's visit to the Bernoulli Center of the \'Ecole
Polytechnique F\'ed\'erale de Lausanne, Switzerland, June-August
2006. The author wishes to express his gratitude to the Center and
its director, professor Tudor Ratiu, for invitation and support.
The author also wants to thank professor Ratiu for the references
\cite{{BR},{DS}} and the discussions concerning their content.
\hspace*{7.5cm}{\small \begin{tabular}{l} Department of
Mathematics\\ University of Haifa, Israel\\ E-mail:
vaisman@math.haifa.ac.il \end{tabular}}
\end{document}